\documentclass[12pt]{article}
\usepackage{amsmath,amsfonts,amsthm,amssymb}
\usepackage{latexsym,amsmath}
\usepackage{graphicx,psfrag}

\numberwithin{equation}{section}

 \DeclareMathOperator{\Ker}{Ker}

 \newcommand{\C}{\mathbb{C}}
 
 \newcommand{\Q}{\mathbb{Q}}
 
 \newcommand{\Z}{\mathbb{Z}}
 \newcommand{\B}{\mathcal{B}}
 \newcommand{\U}{\mathcal{U}}

 \newtheorem{thm}{Theorem}

 \newtheorem{cor}[thm]{Corollary}

 \theoremstyle{definition}

 \newtheorem{rem}{Remark}
 \theoremstyle{remark}
 \newtheorem{clm}{Claim}

\begin{document}

\begin{center}
{\huge No embeddings of solenoids into  surfaces}
\end{center}

\vskip 1 truecm
\begin{center}
BOJU JIANG, jiangbj@math.pku.edu.cn
\\
{\it Depart. of Mathematics, Peking University, Beijing  100871,
China\\}

\end{center}

\begin{center}
SHICHENG WANG, wangsc@math.pku.edu.cn\\
{\it Depart. of Mathematics, Peking University, Beijing 100871,
China\\}
\end{center}

\begin{center}
HAO ZHENG, zhenghao@sysu.edu.cn\\
{\it Depart. of Mathematics, Zhongshan University, Guangzhou
510275, China\\}
\end{center}

{\abstract {\it A quick proof of Bing's theorem indicated by the
title is given. The proof also concludes Gumerov's result on
covering degrees of solenoids.}}

\vskip 1 truecm
 In this paper we always assume that a solenoid is not the
circle.

In his two important papers  on solenoids, Bing proved first that
no solenoid is planar [B1] and then if a solenoid can be embedded
into surfaces, then it must be planar [B2], therefore no solenoid
can be embedded into surfaces. We will give a short proof of this
result. The proof also concludes a recent result of [Gu] on
covering degrees of solenoids.

Identify $S^1$ with the abelian Lie group $U_1 = \{ z\in\C \mid
|z|=1 \}$ and let $\phi_n : S^1 \to S^1$ be the homomorphism defined
by $x \mapsto x^{w_n}$ where $w_n>1$ is an integer. For the sequence
$\{\phi_n: S^1\to S^1\}_{n\ge1}$, its inverse limit, which is
defined as the subspace
$$\Theta = \{(x_0,x_1,...,x_n,...) \mid x_n\in S^1, \; x_{n-1}=\phi_n(x_n)\}$$
of the product space $\Pi_{n=0}^\infty S^1$, is called the {\it
solenoid} of type $\varpi=(w_1,w_2,...,w_n,...)$.  By definition
$\Theta$ is a connected, close and hence compact subgroup of the
abelian topological group $\Pi_{n=0}^\infty S^1$. For more details
see [Mc].

\begin{thm}\label{thm:non-surface}
No solenoid $\Theta$ can be  embedded into a surface.
\end{thm}

\begin{thm}\label{thm:covering}
If $X$ is a finite-fold covering of a solenoid $\Theta$ then each
component of $X$ is homeomorphic to $\Theta$.
\end{thm}

\begin{rem}
According to the description in [Gu, p.2775], seems Theorem
\ref{thm:covering} can be derived from several papers of Fox, Moore
in 1970's, and of Grigorian, Gumerov recently. We will give a direct
proof of Theorem \ref{thm:covering}.
\end{rem}

\begin{proof}[Proof of Theorem \ref{thm:non-surface}]
Suppose there is an embedding $\Theta\subset F$ of a solenoid
$\Theta$ into a surface $F$. We have to seek a contradiction.

First, we may assume $F$ is a closed surface. If not so, since
$\Theta$ is compact, $\Theta$ is contained in the interior of some
connected, compact subsurface $F'\subset F$. Capping a disc on each
component of $\partial F'$, we get an embedding of $\Theta$ into a
closed surface.

Further, we may assume $F$ is orientable. If $F$ is not orientable,
we consider the orientable double covering $\pi: \tilde{F} \to F$.
Since $\pi^{-1}(\Theta)$ is a double covering of $\Theta$, by
Theorem \ref{thm:covering} each component of $\pi^{-1}(\Theta)$ is
homeomorphic to $\Theta$ hence we get an embedding $\Theta \subset
\tilde{F}$ into a closed, orientable surface.

Suppose $\Theta$ is of type $\varpi=(w_1,w_2,...,w_n,...)$. It is
known that the \v{C}ech cohomology group $\check{H}^1(\Theta;R)$
is the direct limit of $R \overset {w_1}\longrightarrow R \overset
{w_2}\longrightarrow R\overset{w_3}\longrightarrow
R\overset{w_4}\longrightarrow \cdots$ (cf. [Mu, p.444]) for any
coefficient ring $R$. It follows that $\check{H}^1(\Theta;\Z)$ is
an infinitely generated $\Z$-module but $\check{H}^1(\Theta;\Q)
\cong \Q$ is a finitely generated $\Q$-module.

By the Alexander duality $\check{H}^p(\Theta;R) \cong
H_{2-p}(F,F\setminus\Theta;R)$ (cf. [GH, p.233]), we have the
exact sequence
$$\cdots \to H_1(F;R) \to \check{H}^1(\Theta;R) \to H_0(F\setminus\Theta;R) \to H_0(F;R) \to \cdots.$$
Since $H_*(F;R)$ is a finitely generated $R$-module for any
commutative ring $R$, the above exact sequence implies that
$H_0(F\setminus\Theta;\Z)$ is an infinitely generated $\Z$-module
and $H_0(F\setminus\Theta;\Q)$ is a finitely generated
$\Q$-module. But on the other hand $H_0(F\setminus\Theta;\Z)$ is a
free abelian group and $H_0(F\setminus\Theta;\Q) =
H_0(F\setminus\Theta;\Z) \otimes_\Z \Q$, so
$H_0(F\setminus\Theta;\Z)$ and $H_0(F\setminus\Theta;\Q)$ have the
same generating set, a contradiction.
\end{proof}

\begin{proof}[Proof of Theorem \ref{thm:covering}] Suppose
$\Theta$ is a solenoid of type $\varpi=(w_1,w_2,...,w_n,...)$ and
$\pi: X \to \Theta$ is an $r$-fold covering.

Let $\Gamma_n$ be the kernel of the projection $p_n: \Theta \to
S^1$ onto the $n$th coordinate, i.e.
$$\Gamma_n = \{ (x_0,x_1,\dots,x_n,\dots) \in \Theta \mid x_0=x_1=\cdots=x_n=1 \}.$$
We have the infinite sequence of close, and hence compact, subgroups
(indeed each $\Gamma_n$ is homeomorphic to the Cantor set, cf. [Mc])
$$\Theta > \Gamma_0 > \Gamma_1 > \Gamma_2 > \cdots > \Gamma_n > \cdots.$$

The right multiplication of $\Gamma_n$ on $\Theta$ makes the
projection $p_n: \Theta \to S^1$ a principal fiber bundle (see also
[Mc]). It follows that $\Theta$ is the mapping tours of a left
transformation $\psi_n: \Gamma_n \to \Gamma_n$. We choose the left
transformation $\psi_n$ as follows. Fix a closed path $\alpha(t) =
e^{2\pi it}$, $t \in [0,1]$, in $S^1$ then $\alpha_n(t) = \gamma(w_1
\cdots w_nt)$ is a lift of $\alpha$ via the projection $p_n$, where
$$\gamma(t) = (e^{2\pi it}, e^{2\pi it/w_1}, e^{2\pi it/w_1w_2}, \dots, e^{2\pi it/w_1 \cdots w_n}, \dots)$$
is the one-parameter subgroup of $\Theta$. Let $\psi_n$ be given
by $x \mapsto \alpha_n(1)x$.

Moreover, as the covering space of $\Theta$, $X$ is the mapping
torus of a lift $\tilde\psi_n$ of $\psi_n$ which is uniquely
determined as follows. For each $\tilde{x} \in \pi^{-1}(\Gamma_n)$,
let $\tilde\alpha$ be the unique lift starting from $\tilde{x}$ of
the path $\alpha_n(t)\pi(\tilde{x})$ in $\Theta$ via the covering
map $\pi$, then $\tilde\psi_n(\tilde{x}) = \tilde\alpha(1)$. Note
that $\psi_n = \psi_0^{w_1 \cdots w_n}|_{\Gamma_n}$ and hence
$\tilde\psi_n = \tilde\psi_0^{w_1 \cdots
w_n}|_{\pi^{-1}(\Gamma_n)}$.

\begin{clm}
Let $\U$ be an open covering of $\Gamma_0$. Then for sufficient
large $n$ each coset of $\Gamma_n$ in $\Gamma_0$ is contained in
some element of $\U$.
\end{clm}

\begin{proof}
Note that $\B = \{ x\Gamma_n \mid n\ge0, \; x\in\Gamma \}$ forms a
basis for the open sets of $\Gamma_0$. So the covering $\U$ has a
refinement $\U'$ which consists of elements of $\B$. Since
$\Gamma_0$ is compact we may assume $\U'$ is finite. Therefore,
for sufficient large $n$ each coset of $\Gamma_n$ in $\Gamma_0$ is
contained in some element of $\U'$ hence in some element of $\U$.
\end{proof}

\begin{clm}
We have $\pi^{-1}(\Gamma_0) \cong \Gamma_0 \times \{1,2,\dots,r\}$
as covering spaces of $\Gamma_0$.
\end{clm}

\begin{proof}
Let $\U$ be the open covering of $\Gamma_0$ which consists of the
open fundamental regions of the covering map
$\pi|_{\pi^{-1}(\Gamma_0)}$ and let $n$ be sufficiently large as
in Claim 1. Then each coset $x\Gamma_n$ in $\Gamma_0$ is contained
in some element $\U$ hence $\pi^{-1}(x\Gamma_n) \cong x\Gamma_n
\times \{1,2,\dots,r\}$ as covering spaces of $x\Gamma_n$. Since
the cosets of $\Gamma_n$ are disjoint and open in $\Gamma_0$, the
claim follows.
\end{proof}

In what follows we fix a homeomorphism from Claim 2 and identify
both sets.

\begin{clm}
 For each sufficient large $n$ there exists a permutation
$\sigma_n$ of $\{1,2,\dots,r\}$ such that $\tilde\psi_n(x,j) =
(\psi_n(x),\sigma_n(j))$.
\end{clm}

\begin{proof}
Applying Claim 1 on the the open covering of $\Gamma_0$ which
consists of the open (and close) sets $U_\sigma = \{ x\in\Gamma_0
\mid \tilde\psi_0(x,j) = (\psi_0(x),\sigma(j)) \}$ with $\sigma$
running over all permutations of $\{1,2,\dots,r\}$, one notices
that for sufficiently large $n$ each coset of $\Gamma_n$ in
$\Gamma_0$ is contained in some $U_\sigma$. Since $\tilde\psi_n =
\tilde\psi_0^{w_1 \cdots w_n}|_{\pi^{-1}(\Gamma_n)}$, the claim
follows.
\end{proof}

By Claim 3 $X$ is the disjoint union of the mapping tori of
$\tilde\psi_n|_{\Gamma_n \times J}$ where $J$ runs over all
$\sigma_n$-orbits. Since $\Theta$ is connected, the following two
claims eventually establish the theorem. Below we denote by $T_f$
the mapping torus of a self homeomorphism $f$.

\begin{clm}
For sufficient large $n$ the length of each $\sigma_n$-orbit is
relatively prime to $w_{n'}$ for all $n'>n$.
\end{clm}

\begin{proof}
Note that $\sigma_{n+1} = \sigma_n^{w_{n+1}}$. Therefore, each
$\sigma_{n'}$-orbit is contained in some $\sigma_n$-orbit for
$n'>n$ and if the length of a $\sigma_n$-orbit $J$ is not
relatively prime to $w_{n'}$ for some $n'>n$ then $J$ splits into
several $\sigma_{n'}$-orbits. Since a permutation of
$\{1,2,\dots,r\}$ has at most $r$ orbits, the claim follows.
\end{proof}

\begin{clm}
Let $n$ be sufficient large as in Claim 3. If the length of a
$\sigma_n$-orbit $J$ is relatively prime to $w_{n'}$ for all
$n'>n$, then $T_{\tilde\psi_n|_{\Gamma_n \times J}} \cong \Theta$.
\end{clm}

\begin{proof}
Let $l = |J|$ be the length of $J$. It is clear that
$T_{\tilde\psi_n|_{\Gamma_n \times J}} \cong T_{\psi_n^l}$. Note
that $\Gamma_n$ is the inverse limit of the sequence
$$\{ \phi_{n+k}: \Ker(\phi_{n+1}\cdots\phi_{n+k}) \to \Ker(\phi_{n+1}\cdots\phi_{n+k-1}) \}_{k\ge1}.$$
Since $l$ is relatively prime to $w_{n'}$ for all $n'>n$, the
homomorphisms
$$\Ker(\phi_{n+1}\cdots\phi_{n+k}) \to \Ker(\phi_{n+1}\cdots\phi_{n+k})$$
defined by $x \mapsto x^l$ are isomorphic. It follows that the
homomorphism $\Gamma_n \to \Gamma_n$ defined by $x \to x^l$ is
isomorphic, via which one notices that $\psi_n$ is topologically
conjugate to $\psi_n^l$. So we have $T_{\psi_n} \cong T_{\psi_n^l}$
and therefore, $T_{\tilde\psi_n|_{\Gamma_n \times J}} \cong
T_{\psi_n^l} \cong T_{\psi_n} \cong \Theta$.
\end{proof}

From Claim 4 and Claim 5 we also have

\begin{cor}[{[Gu]}]\label{cor:existence}
A solenoid of type $\varpi=(w_1,w_2,...,w_n,...)$ has a connected
$r$-fold covering if and only if $r$ is relatively prime to all
but finitely many $w_n$.
\end{cor}
\end{proof}

\bibliographystyle{amsalpha}

\end{document}